\documentclass[12pt]{siamart220329}

\usepackage{tim,epsfig,fullpage,cite,algorithm,algorithmic}
\usepackage{ma580}
\usepackage{graphicx}
\usepackage{hyperref}
\usepackage{xcolor,mdframed}
\usepackage{url}
\def\calf{{\cal F}}

\newcommand\audiencename{Audience}

\title{Interprecision transfers in iterative refinement}
\author{
C. T. Kelley%
\thanks{North Carolina State University,
Department of Mathematics,
Box 8205, Raleigh, NC 27695-8205, USA
(Tim\_Kelley@ncsu.edu).
This work was partially supported by
the Center for Exascale Monte-Carlo Neutron Transport (CEMeNT) a 
PSAAP-III project funded by Department of Energy
grant number  DE-NA003967.
}
}

\begin{document}\maketitle

\begin{abstract}
We make the interprecision transfers explicit in an algorithmic
description of iterative refinement and obtain new insights into
the algorithm. One example is 
the classic variant of iterative refinement where the
matrix and the factorization are stored in a working precision and
the residual is evaluated in a higher precision. In that case
we make the observation
that this algorithm will solve a promoted form of the original problem
and thereby characterize the limiting behavior in a novel way and obtain
a different version of the classic convergence
analysis. We also discuss two approaches
for interprecision transfer in the triangular solves.
\end{abstract}

\begin{keywords}
Iterative refinement, Interprecision transfers,
Mixed-precision arithmetic, Linear systems
\end{keywords}

\begin{AMS}
65F05, 
65F10, 
\end{AMS}


\section{Introduction}
\label{sec:intro}

Iterative refinement (IR) is a way to lower factorization costs
in the numerical solution of a linear system $\ma \vx = \vb$
by performing the factorization in a lower precision. 
Algorithm~{\bf IR-V0} 
is a simple formulation using Gaussian elimination. In this
formulation all computations are done in a high precision
except for the $LU$ factorization. 

\begin{algorithm}
{$\mbox{\bf IR-V0}(\ma, \vb)$}
\begin{algorithmic}
\STATE $\vx = 0$
\STATE $\vr = \vb$
\STATE Factor $\ma = \ml \mU$ in a lower precision
\WHILE{$\| \vr \|$ too large}
\STATE $\vd = \mU^{-1} \ml^{-1} \vr$
\STATE $\vx \leftarrow \vx + \vd$
\STATE $\vr = \vb - \ma \vx$
\ENDWHILE
\end{algorithmic}
\end{algorithm}

Algorithm~{\bf IR-V0} leaves out many implementation details. Some
recent papers 
\cite{amestoy:2024,CarsonHigham1,CarsonHigham,demmelir} have made the
algorithmic details explicit. The purpose of this paper is to build
upon that work by explicitly including the interprecision transfers
in the algorithmic description. One consequence of this, which
we discuss in \S~\ref{subsec:resid} and 
\S~\ref{subsec:residhigh} is a novel
interpretation of the classic form of the method \cite{Wilkinson48}
where the residual is evaluated in an extended precision.

\subsection{Notation}
\label{subsec:notation}

We use the terminology from \cite{demmelir, amestoy:2024}
and consider several precisions. We will use Julia-like notation for
data types.

\begin{itemize}
\item The matrix $\ma$ and right side $\vb$ are stored in the
{\bf working precision} $TW$.
\item $\ma$ is factored in the {\bf factorization precision} $TF$.
\item The residual is computed in the {\bf residual precision} $TR$.
\item The triangular solves for $\ml \mU \vd = \vr$ are done in the
{\bf solver precision} $TS$.
\end{itemize}

If TH is a higher precision that TL we will write $TL < TH$. $fl_X$
will be the rounding operation to precision $TX$. We will let
$\calf_{X}$ be the floating point numbers with precision $TX$.
We will assume that
a low precision number can be exactly represented in any higher
precision. So
\begeq
\label{eq:nestedpre}
x \in \calf_X \subset \calf_Y
\mbox{ if $TX \le TY$}.
\endeq

We will let $u_X$ denote the unit roundoff for precision $TX$ and let
$I_{A}^{B}$ denote the interprecision transfer from precision $TA$ to
precision $TB$. Interprecision transfer is more than rounding and can,
in some cases, include data allocation. If $TA < TB$ ($TA > TB$)
then we will call the transfer $I_A^B$ {\em upcasting}
({\em downcasting}). 

Upcasting is simpler than downcasting because if $TA < TB$ and 
$x \in \calf_{A}$, then \eqnok{nestedpre} implies that
\begeq
\label{eq:isoxf}
I_A^B (x) = x.
\endeq
However, upcasting is not linear. In fact if $y \in \calf_{A}$ there is
no reason to expect that
\[
I_A^B (xy) = I_A^B(x) I_A^B(y)
\]
because the multiplication on the left is done in a lower precision than
the one on the right. 
Downcasting is more subtle because not only is it nonlinear but also
if $TA  > TB$ and $x \in \calf_A$,
then \eqnok{isoxf} holds only if $x \in \calf_B \subset \calf_A$. 

The nonlinearity of interprecision transfers should be made explicit in
analysis especially for the triangular solves (see \S~\ref{subsec:solve}).

\section{Interprecision transfers in IR}
\label{sec:results}
Algorithm~{\bf IR0} includes many implicit interprecision transfers. We
will make the assumption, which holds in IEEE \cite{IEEEnew} arithmetic,
that when a binary operation $\circ$ is performed between floating
point numbers with different precisions, then the lower precision number
is promoted before the operation is performed.

So, if $TH > TL$, $u \in \calf_{H}$, $v \in \calf_{L}$, then
\begeq
\label{eq:rule}
fl_{H} (u \circ v) = fl_{H} (u \circ (I_{L}^{H} v))
\mbox{ and }
fl_{H} (v \circ u) = fl_{H} ((I_{L}^{H} v) \circ u).
\endeq

We will use \eqnok{rule} throughout this paper when we need to
make the implicit interprecision transfers explicit.

\subsection{The low precision factorization}
\label{subsec:factor}

We will begin with the low precision factorization. The
line in Algorithm~{\bf IR-V0} for this is

\begin{itemize}
\item Factor $\ma = \ml \mU$ in a lower precision.
\end{itemize}

However, $\ma$ is stored in precision $TW$
and the factorization is in
precision $TF$. Hence one must make a copy of $\ma$. So to make this
interprecision transfers explicit we should express this as

\begin{itemize}
\item Make a low precision copy of $\ma$. $\ma_F = I_{W}^{F} \ma$.
\item Compute an $LU$ factorization $\ml \mU$ of $\ma_F$, overwriting
$\ma_F$.
\end{itemize}

So in this way the storage costs of IR become clear and we can
see the time/storage tradeoff. If, for example, $TW = Float64$ and
$TF = Float32$, one must allocate storage for $\ma_F$ which is a
50\%  increase in matrix storage.

\subsection{The residual precision}
\label{subsec:resid}

For the remainder of this paper we use the $\ell^\infty$ norm, so
$\| \cdot \| = \| \cdot \|_\infty$.

The algorithm {\bf IR-V0} does not make it clear how the residual precision
affects the iteration. The description in \cite{demmelir} carefully
explains what one must do. 

\begin{itemize}
\item Store $\vx$ in precision $TR$.
\item Solve the triangular system in precision $TS=TR$.
\item Store $\vd$ in precision $TR$.
\end{itemize}

The most important consequence of this and \eqnok{rule} is
that the residual is
\[
\vr = I_W^R (\vb - \ma \vx) = (I_W^R \ma) \vx - I_W^R \vb.
\]
While $I_W^R \ma = \ma$ by our assumptions that $TR \ge TW$, the 
residual is computed in precision $TR$ and is therefore the residual
of a promoted problem and the iteration is approximating 
the solution of that promoted problem
\begeq
\label{eq:promotesol}
\vx_P^* = (I_W^R \ma)^{-1} \vb = (I_W^R \ma)^{-1} I_W^R \vb,
\endeq
which is posed in precision TR. We make a distinction between
$\vx_P^*$ and $\vx^* = \ma^{-1} \vb$ only in those cases, such as
\eqnok{promotesol} where we are talking about the computed solution
in the residual precision. In exact arithmetic, of course,
$\vx_P^* = \vx^*$.

All of these interprecision transfers are implicit and need not be done
within the iteration. For example, when one computes $\vr$ in precision
$TR$, then \eqnok{rule} implies that
the matrix-vector product $\ma \vx$ automatically promotes the
elements of $\ma$ to precision $TR$ because $\vx$ is stored in precision
TR.

We can use the fact that convergence is to the solution of the 
promoted problem to get a simple error estimate for
the classic special case \cite{Wilkinson48}. Here
$TR = TS > TW \ge TF$. If the iteration terminates with
\[
\| \vr \|/\| \vb \| \le \tau
\]
then the standard estimates \cite{demmel} imply that
\begeq
\label{eq:showwilk}
\frac{\| \vx - \vx^* \|}{\| \vx^*\|}
= \frac{\| \vx - \vx^*_P \|}{\| \vx^*_P \|}
\le \frac{ \kappa(I_W^R \ma) \| \vr \|}{\| I_W^R (\vb) \|}
= \frac{ \kappa(\ma) \| \vr \|}{\| \vb \|}
\le \tau \kappa(\ma).
\endeq
In \eqnok{showwilk}
we use the fact that $\| I_W^R (\vb) \| = \| \vb \|$ in the 
$\ell^\infty$ norm and use the exact value of $\kappa(\ma)$, which
is the same as $\kappa(I_W^R \ma)$ by \eqnok{nestedpre}.
The estimate \eqnok{showwilk} is also
true if $TR = TW$, but then there is no need to consider a promoted 
problem.

\subsection{The triangular solve}
\label{subsec:solve}

The choice of $TS$ affects the number of interprecision transfers 
and the storage cost in the
triangular solve. The line in Algorithm~{\bf IR-V0} for this is 

\begin{itemize}
\item $\vd = \mU^{-1} \ml^{-1} \vr$
\end{itemize}

If $TS = TR$, then the LAPACK's default behavior
is to do the interprecision transfers as needed using \eqnok{rule}. The
subtle consequence of this is that 
\[
(\ml \mU)^{-1} \vr = ( (I_{F}^{R} \ml) (I_{F}^{R} \mU ) )^{-1} \vr.
\]
Hence, one is implicitly doing the triangular solves with the factors
promoted to the residual precision. We will refer to this approach
as {\bf on-the-fly} interprecision transfers.

Combining the results with \S~\ref{subsec:resid} we can expose all the
interprecision transfers in the transition
from a current iteration $\vx_c$ to a new one $\vx_+$. In the case
$TS = TR$ we have a linear stationary iterative method.
The computation is done entirely in $TR$.
\begeq
\label{eq:itconv}
\vx_+  = \vx_c + \vd = \vx_c + (\ml \mU)^{-1} 
(I_W^R \vb - \ma \vx_c) = \mm \vx_c +  (\ml \mU)^{-1} I_W^R \vb,
\endeq
where the iteration matrix is
\begeq
\label{eq:itmat}
\mm = \mi - ((I_{F}^{R} \ml) (I_{F}^{R} \mU) )^{-1} (I_W^R \ma).
\endeq

The residual update is
\begeq
\label{eq:resconv}
\vr_+  = \mm_r \vr_c
\endeq
where
\begeq
\label{eq:itmatr}
\mm_r = \mi -  (I_W^R \ma) ((I_{F}^{R} \ml) (I_{F}^{R} \mU) )^{-1}.
\endeq
One must remember that if $TR > TW$ and $\vx \in \calf_W^N$
then $(I_W^R \ma) \vx \ne \ma \vx$ because
the matrices are in different precisions and matrix-vector products produce
different results.

All the interprecision transfers
in \eqnok{itmat} and \eqnok{itmatr}
are implicit and the promoted matrices
are not actually stored. However, the
promotions matter because they can help avoid underflows and overflows
and influence the limit of the iteration.

If $TS = TF < TW$, then the interprecision transfer is done before
the triangular solves and the number of interprecision transfers is
$N$ rather than $N^2$. We will refer to this as interprecision transfer
{\bf in-place} to distinguish it from on-the-fly.
For in-place interprecision transfers  we copy $\vr$ from the residual
precision $TR$ to the factorization precision before the solve and
then upcast
the output of the solve back to the residual precision. 
So, one must store the low precision copy of $\vr$
In this case one should scale $\vr$ before the downcasting transfer
$I_R^F$ \cite{highamscaling}. One reason for this is that the 
absolute size of $\vr$ could be very small, as would be the case in the
terminal phase of IR, and one could underflow before the iteration is
complete. 
So the iteration in this case is
\begeq
\label{eq:ipir}
\vx_+ = \vx_c + \vd = \vx_c + \| \vr \|  I_F^R \left[ (\ml \mU)^{-1} 
\frac{I_R^F \vr}{\| \vr \|} \right].
\endeq
This is not a stationary linear iterative method because the map
$\vr \rightarrow \frac{I_R^F \vr}{\| \vr \|}$ is nonlinear.

Even though downcasting $\vr$ reduces the interprecision transfer cost,
one should do interprecision transfers on-the-fly if $TF$ is half precision,
if $TR > TW$, or if one is using the low-precision factorization as a
preconditioner \cite{amestoy:2024,CarsonHigham}.

\section{Explicit Interprecision Transfers}
\label{sec:explicit}

We apply the results from \S~\ref{sec:results} 
Algorithm~{\bf IR-V0} and obtain Algorithm~{\bf IR-V1}, where
all the interprecision transfers are explicit. 

\begin{algorithm}
{$\mbox{\bf IR-V1}(\ma, \vb, TF, TW, TR)$}
\begin{algorithmic}
\STATE $\vx = 0 \in \calf_R^N$
\STATE $\vr = I_W^R(\vb)$
\STATE $\ma_F = I_W^F( \ma)$.
\STATE Factor $\ma_F = \ml \mU$ in precision TF
\WHILE{$\| \vr \|$ too large}
\STATE $\vd = ((I_F^R \ml) (I_F^R \mU) )^{-1} \vr$ in precision TR
\STATE $\vx \leftarrow \vx + \vd$ 
\STATE $\vr = (I_W^R \vb) - (I_W^R \ma) \vx$
\ENDWHILE
\end{algorithmic}
\end{algorithm}

In the remainder of this section we look at some consequences of
this formulation of IR.
 
\subsection{The case $TS = TR > TW$}
\label{subsec:residhigh}

In this section we will assume that the triangular solves
are done in the residual precision (TS = TR). In this case, one can
see from Algorithm~{\bf IR-V1} that no computations are done in the
working precision at all. The working precision is only used to store
$\ma$ and $\vb$, but residual computations are done in the residual precision
with promotion on-the-fly. We state this observation as a theorem.

\begin{theorem}
\label{th:residhigh}
If $TS=TR > TW$, then the three precision algorithm
\[
\mbox{\bf IR-V1}(\ma, \vb, TF, TW, TR)
\]
produces the same computed results as the two precision algorithm
\[
\mbox{\bf IR-V1}(I_W^R \ma, I_W^R \vb, TF, TR, TR).
\]
\end{theorem}

The theorem makes it clear that the iteration is reducing the residual of the
promoted problem. So we can apply the classical ideas for {\bf IR-V0}
\cite{higham,higham97}
and understand the case $TR > TW$ in that way. For example, 
if $\ma$ is not highly
ill-conditioned, the LU factorization is stable, and the norm of
iteration matrix for IR $\| \mm \| < 1$, then
we can use equation (4.9) from \cite{higham97}
to obtain
\begeq
\label{eq:linearconv}
\| \vr_+ \| \le G \| \vr_c \| + g, \mbox{ where }
g = O(u_R [\| I_W^R \ma \| \| \vx^*_P \| + \| I_W^R \vb \|] ). 
\endeq
where $G < 1$.  Hence we will be able to 
reduce $\| \vr \|$ until the iteration saturates with 
$\| \vr \| \approx g$. 

Since one has no {\em a priori} knowledge of
$g$, one must manage the iteration in a way to detect stagnation. 
The recommendation from
\cite{higham97} is to terminate the iteration when

\begin{enumerate}
\item $\| \vr \|  \le u_R (\| \ma \| \| \vx \| +  \| \vb \|)$,
\item $\| \vr_+ \| \ge \alpha \| \vr_c \|$, or
\item too many iterations have been performed.
\end{enumerate}

The first item in the list is successful convergence where we approximate
the norms of the promoted objects $I_W^R \ma$ and $I_W^R \vb$ with the 
ones in the working precision which have stored.
The two failure modes are insufficient decrease in the residual 
and slow convergence in the terminal phase.
The recommendation in \cite{higham97} was to set
$C=1$ and $\alpha=.5$, and to limit IR to five iterations. 
Our Julia code \cite{ctk:mparrays} uses a variation of this approach.
We use $\alpha = .9$ in our solver and do not
put a limit on the iterations. The reason for these choices are to give
the IR iteration a better chance to terminate successfully.

So, if we couple the termination strategy with \eqnok{showwilk} we see
that if the iteration terminates successfully then
\[
\frac{\| \vx - \vx^*_P \|}{\| \vx^*_P \|}
\le u_R
\frac{ C (\| \ma \| \| \vx \| +  \| \vb \|) \kappa(\ma) \| \vb \|}%
{\| I_W^R (\vb) \|}.
\]

When $TR = TS > TW$ one can also attempt to estimate the 
convergence rate $G$ from \eqnok{linearconv} and then estimate
the error $\| \vx - \vx^*_P \|$. This is a common strategy in the 
nonlinear solver literature, especially for stiff initial
value problems \cite{slc,vode,daspk,dassl,ctk:roots}.
The idea is that as the iteration progresses
\[
G \approx \sigma = \| \vx_{n+1} - \vx_n \| / \| \vx_n - \vx_{n-1} \|
\]
is a very good estimate if $G$ is small enough. In that case
\[
\| \vx_n - \vx^*_P \| \le \frac{\| \vx_{n+1} - \vx_n \|}{1 - \sigma} 
\]
and one can terminate the iteration when one predicts that
\[
\| \vx_{n+1} - \vx^*_P \| \le \frac{\| \vx_{n+1} -  \vx_n \| \sigma}{1 - \sigma}
\]
is sufficiently small. The algorithm in \cite{demmelir} 
does this and terminates when the predicted error is less than $u_W$.

\subsection{Cost of Interprecision Transfers}
\label{subsec:ircost}

One case where setting $TS = TF$ may be useful is if $TS=Float32$
and $TW = TR = Float64$. In this case, unlike $TF = Float16$, 
tools such as LAPACK and BLAS have been compiled to work efficiently. 
The cases of interest in this section are medium sized
problems where the $O(N^3)$ cost of the LU factorization is a few times
more than the cost of the triangular solves, but not orders of magnitude more.
In these cases the $O(N^2)$ cost of interprecision transfers on the fly
is noticeable and could make a difference in cases where many triangular
solves are done for each factorization. One example is for nonlinear
solvers \cite{ctk:roots,ctk:fajulia} where the factorization of the
Jacobian can be reused for many Newton iterations (or even time
steps when solving stiff initial value problems \cite{dassl,shampinebook}).

We illustrate this with some cpu timings. The computations in this section
were done on an Apple Macintosh 
Mini Pro with a M2 processor and eight performance
cores. We used OpenBlas, which satisfies \eqnok{rule}, rather than the
AppleAccelerate Framework, which does not. We used
Julia \cite{Juliasirev} v1.11.0-beta2 with the author's 
{\bf MulitPrecisionArrays.jl} \cite{ctk:mparrays,ctk:mparraysdocs}
Julia package. We made this choice because
Julia v1.11.0 has faster matrix-vector products
than the current version v1.10.4.

We used the Julia package {\bf BenchmarkTools.jl}
\cite{benchmarktools} to get the timings we report in 
Table~\ref{tab:otf}.  This is the standard way to obtain timings in Julia.
BenchmarkTools repeats computations and can obtain accurate results even if
the compute time per run is very small.

We have put the Julia codes that generate Table~\ref{tab:otf} 
from \S~\ref{subsubsec:example} in a
GitHub repository

\url{https://github.com/ctkelley/IR_Precision_Transfers}

\subsubsection{Integral Equation Example}
\label{subsubsec:example}

We will use a concrete example rather than generating random problems. For a
given dimension $N$ let $\mg$ the matrix corresponding to the composite
trapezoid rule discretization of the Greens operator $\calg$
for $-d^2/dx^2$ on $[0,1]$

\[
\calg u(x) = \int_0^1 g(x,y) u(y) \, dy \mbox{ where }
g(x,y) =
    \left\{\begin{array}{c}
        y (1-x) ; \ x > y\\
        x (1-y) ; \ x \le y
    \end{array}\right.
\]

The eigenvalues of $\calg$ are $1/(n^2 \pi^2)$ for $n = 1, 2, \dots$. 

We use $\ma = \mi - 800.0 * \mg$ in this example. The conditioning of
$\ma$ is somewhat poor with an $\ell^\infty$ condition number of roughly
$\kappa_\infty (\ma) \approx 18,253$.
for the dimensions we consider in this section.

We terminate the IR iteration with the residual condition from 
\S~\ref{subsec:residhigh} and tabulate the dimension, the time (LU) for
copying $\ma$ from TW to TF and performing the 
LU factorization in TF (column 2 of Table~\ref{tab:otf},
the timings for the two variants of the 
triangular solves (OTF = on-the-fly: column 3, IP = in-place: column 4),
and the times
and iteration counts for the IR loop with the two variations of the
triangular solves (columns 5--8).

\begin{table}[h]
\caption{\label{tab:otf} Cost of OTF Triangular Solve}
\centerline{
\begin{tabular}{llllllll} 
        N &       LU &      OTF &       IP &   OTF-IR &      its &    IP-IR &      its \\ 
\hline 
200 & 1.6e-04 & 1.5e-05 & 5.7e-06 & 4.5e-05 & 3 & 3.8e-05 & 3   \\ 
400 & 4.9e-04 & 5.4e-05 & 1.9e-05 & 2.3e-04 & 4 & 2.3e-04 & 5   \\ 
800 & 1.8e-03 & 2.4e-04 & 6.6e-05 & 9.9e-04 & 5 & 7.3e-04 & 5   \\ 
1600 & 7.8e-03 & 1.4e-03 & 2.5e-04 & 2.8e-03 & 4 & 2.1e-03 & 4   \\ 
3200 & 4.2e-02 & 9.2e-03 & 1.3e-03 & 2.0e-02 & 5 & 1.7e-02 & 5   \\ 
6400 & 2.9e-01 & 3.6e-02 & 6.1e-03 & 7.2e-02 & 5 & 5.8e-02 & 5   \\ 
\hline 
\end{tabular} 
}
\end{table}

Table~\ref{tab:otf} shows that the
factorization time is between 6 and 10 times that of the on-the-fly triangular
solve indicating that the triangular solves could be significant
if the matrix-vector products were fast or one had to solve for many 
right hand sides.  We saw that effect in the nonlinear
examples in \cite{ctk:sirev20,ctk:fajulia} where the nonlinear residual
could be evaluated in $O(N log N)$ work.
One can also see that
the cpu time for the in-place
triangular solves is 2--5 times less than for the on-the-fly version. 

In the final four columns we see that,
while the version with in-place triangular solves is somewhat faster, the
difference is not compelling. This is no surprise because the 
matrix-vector product in the residual precision (double) takes $O(N^2)$ work
and is therefore a 
significant part of the cost for each IR iteration.  The number of
iterations is the same in all but one case, with in-place triangular
solves taking one more iteration in that case. That is consistent with
the prediction in \cite{CarsonHigham}.

\section{Conclusions}

We expose the interprecision transfers in iterative refinement and obtain
new insights into this classical algorithm. In particular we show that
the version in which the residual is evaluated in an extended precision is
equivalent to solving a promoted problem and show how interprecision transfers
affect the triangular solves.

\section{Acknowledgments}
The author is very grateful to Ilse Ipsen for listening to him
as he worked through the ideas for this paper.

\bibliographystyle{siamplain}
\bibliography{IR_Note}

\begin{thebibliography}{10}

\bibitem{amestoy:2024}
{\sc P.~Amestoy, A.~Buttari, N.~J. Higham, J.-Y. L’Excellent, T.~Mary, and
  B.~Vieubl\'{e}}, {\em Five-precision gmres-based iterative refinement}, SIAM
  Journal on Matrix Analysis and Applications, 45 (2024), pp.~529--552,
  \url{https://doi.org/10.1137/23M1549079}.

\bibitem{Juliasirev}
{\sc J.~Bezanson, A.~Edelman, S.~Karpinski, and V.~B. Shah}, {\em Julia: A
  fresh approach to numerical computing}, SIAM Review, 59 (2017), pp.~65--98.

\bibitem{slc}
{\sc K.~E. Brenan, S.~L. Campbell, and L.~R. Petzold}, {\em The Numerical
  Solution of Initial Value Problems in Differential-Algebraic Equations},
  no.~14 in Classics in Applied Mathematics, SIAM, Philadelphia, 1996.

\bibitem{vode}
{\sc P.~N. Brown, G.~D. Byrne, and A.~C. Hindmarsh}, {\em {VODE}: A variable
  coefficient ode solver}, SIAM J. Sci. Statist. Comput., 10 (1989),
  pp.~1038--1051.

\bibitem{daspk}
{\sc P.~N. Brown, A.~C. Hindmarsh, and L.~R. Petzold}, {\em Using {K}rylov
  methods in the solution of large-scale differential-algebraic systems}, SIAM
  J. Sci. Comput., 15 (1994), pp.~1467--1488.

\bibitem{CarsonHigham1}
{\sc E.~Carson and N.~J. Higham}, {\em A new analysis of iterative refinement
  and its application of accurate solution of ill-conditioned sparse linear
  systems}, SIAM Journal on Scientific Computing, 39 (2017), pp.~A2834--A2856,
  \url{https://doi.org/10.1137/17M112291}.

\bibitem{CarsonHigham}
{\sc E.~Carson and N.~J. Higham}, {\em Accelerating the solution of linear
  systems by iterative refinement in three precisions}, SIAM Journal on
  Scientific Computing, 40 (2018), pp.~A817--A847,
  \url{https://doi.org/10.1137/17M1140819}.

\bibitem{benchmarktools}
{\sc J.~{Chen} and J.~{Revels}}, {\em {Robust benchmarking in noisy
  environments}}, 2016, \url{https://arxiv.org/abs/1608.04295}.

\bibitem{demmelir}
{\sc J.~Demmel, Y.~Hida, W.~Kahan, X.~S. Li, S.~Mukherjee, and E.~J. Riedy},
  {\em Error bounds from extra-precise iterative refinement}, ACM Trans. Math.
  Soft.,  (2006), pp.~325--351.

\bibitem{demmel}
{\sc J.~W. Demmel}, {\em Applied Numerical Linear Algebra}, SIAM, Philadelphia,
  1997.

\bibitem{higham}
{\sc N.~J. Higham}, {\em Accuracy and Stability of Numerical Algorithms},
  Society for Industrial and Applied Mathematics, Philadelphia, PA, USA, 1996,
  \url{http://www.ma.man.ac.uk/~higham/asna.html}.

\bibitem{higham97}
{\sc N.~J. Higham}, {\em Iterative refinement for linear systems and {LAPACK}},
  IMA J. Numer. Anal., 17 (1997), pp.~495--509.

\bibitem{highamscaling}
{\sc N.~J. Higham, S.~Pranesh, and M.~Zounon}, {\em Squeezing a matrix into
  half precision, with an application to solving linear systems}, SIAM J. Sci.
  Comp., 41 (2019), pp.~A2536--A2551.

\bibitem{IEEEnew}
{\sc {IEEE Computer Society}}, {\em {IEEE} standard for floating-point
  arithmetic, {IEEE} {Std} 754–2019}, July 2019.

\bibitem{ctk:roots}
{\sc C.~T. Kelley}, {\em {Iterative Methods for Linear and Nonlinear
  Equations}}, no.~16 in Frontiers in Applied Mathematics, SIAM, Philadelphia,
  1995.

\bibitem{ctk:sirev20}
{\sc C.~T. Kelley}, {\em Newton's method in mixed precision}, SIAM Review, 64
  (2022), pp.~191--211, \url{https://doi.org/10.1137/20M1342902}.

\bibitem{ctk:fajulia}
{\sc C.~T. Kelley}, {\em {Solving Nonlinear Equations with Iterative Methods:
  Solvers and Examples in Julia}}, no.~20 in Fundamentals of Algorithms, SIAM,
  Philadelphia, 2022.

\bibitem{ctk:mparrays}
{\sc C.~T. Kelley}, {\em {MultiPrecisionArrays.jl}}, 2023,
  \url{https://doi.org/10.5281/zenodo.7521427},
  \url{https://github.com/ctkelley/MultiPrecisionArrays.jl}.
\newblock Julia Package.

\bibitem{ctk:mparraysdocs}
{\sc C.~T. Kelley}, {\em Using {MultiPrecisonArrays.jl}: {I}terative refinement
  in {J}ulia}, 2024, \url{https://arxiv.org/abs/2311.14616}.

\bibitem{dassl}
{\sc L.~R. Petzold}, {\em A description of {DASSL}: a differential/algebraic
  system solver}, in Scientific Computing, {R. S. Stepleman et al.}, ed., North
  Holland, Amsterdam, 1983, pp.~65--68.

\bibitem{shampinebook}
{\sc L.~F. Shampine}, {\em Numerical Solution of Ordinary Differential
  Equations}, Chapman and Hall, New York, 1994.

\bibitem{Wilkinson48}
{\sc J.~H. Wilkinson}, {\em Progress report on the automatic computing engine},
  Tech. Report MA/17/1024, Mathematics Division, Department of Scientific and
  Industrial Research, National Physical Laboratory, Teddington, UK, 1948,
  \url{http://www.alanturing.net/turing_archive/archive/l/l10/l10.php}.

\end{thebibliography}

\end{document}